\documentclass[10pt]{amsart}
\usepackage{amssymb,times,amsmath,cancel,mathrsfs,multirow,textcomp,calligra}
\usepackage{pgfplots,xxcolor}
\usepackage{graphicx,color,epsfig}
\usepackage{bbm}
\usepackage{caption}

\usepackage[all]{xy}
\xyoption{all}
\usepackage{tikz}
\usetikzlibrary{calc,3d,arrows,decorations.markings}
\tikzset{->-/.style={decoration={  markings,  mark=at position .75 with {\arrow{latex}}},postaction={decorate}}}
\tikzset{-<-/.style={decoration={  markings,  mark=at position .75 with {\arrow{latex reversed}}},postaction={decorate}}}

\pgfdeclareradialshading[tikz@ball]{ball}{\pgfqpoint{0bp}{0bp}}{%
 color(0bp)=(tikz@ball!0!white);
 color(7bp)=(tikz@ball!0!white);
 color(15bp)=(tikz@ball!70!black);
 color(20bp)=(black!70);
 color(30bp)=(black!70)}
\makeatother

\tikzset{viewport/.style 2 args={
    x={({cos(-#1)*1cm},{sin(-#1)*sin(#2)*1cm})},
    y={({-sin(-#1)*1cm},{cos(-#1)*sin(#2)*1cm})},
    z={(0,{cos(#2)*1cm})}
}}

\pgfplotsset{only foreground/.style={
    restrict expr to domain={rawx*\CameraX + rawy*\CameraY + rawz*\CameraZ}{-0.05:100},
}}
\pgfplotsset{only background/.style={
    restrict expr to domain={rawx*\CameraX + rawy*\CameraY + rawz*\CameraZ}{-100:0.05}
}}

\def\addFGBGplot[#1]#2;{
    \addplot3[#1,only background, opacity=0.25] #2;
    \addplot3[#1,only foreground] #2;
}

\newcommand{\ViewAzimuth}{-60}
\newcommand{\ViewElevation}{45}

\let\cal\mathcal

\def\cO{{\cal O}}

\def\cS{{\cal S}}

\let\blb\mathbb

\def \PP{{\blb P}}

\def \T{{\blb T}}

\def \MM{{\blb M}}

\def \SSS{{\blb S}}

\newcommand{\weg}[1]{}

\newcommand{\C}{\mathbb{C}}

\newcommand{\R}{\mathbb{R}}

\newcommand{\sM}{\mathscr{M}}
\newcommand{\sO}{\mathscr{O}}

\newcommand{\sm}[1]{\left(\begin{smallmatrix}#1\end{smallmatrix}\right)}

\newcommand{\erbij}[1]{}

\theoremstyle{definition}

{

}

\theoremstyle{remark}

\newcommand{\Ker}{\textrm{Ker}}
\newcommand{\Image }{\textrm{Im}}

\renewcommand{\C}{\mathbb{C}}

\title{Reflections in a cup of coffee}

\author{Raf Bocklandt}
\address{Raf Bocklandt\\
Korteweg de Vries institute\\
University of Amsterdam (UvA)\\
Science Park 904\\ 
1098 XH Amsterdam\\ 
The Netherlands
}
\email{raf.bocklandt@gmail.com}

\xyoption{all}

\begin{document}
\begin{abstract}
Allegedly, Brouwer discovered his famous fixed point theorem while stirring a cup of coffee and noticing that there is always
at least one point in the liquid that does not move. In this paper, based on a talk in honour of Brouwer at the University of Amsterdam, we will explore how Brouwer's ideas about this phenomenon spilt over in a lot
of different areas of mathematics and how this eventually led to an intriguing geometrical theory we now know as mirror symmetry.
\end{abstract}
\maketitle

\section{Brouwer's beginnings}
As the nineteenth century drew to a close,
the famous French mathematician Henri Poincar\'e unleashed a new type of geometry onto the world, which would transform many areas in mathematics and keep 
mathematicians busy until this very day. Poincar\'e used topology ---or analysis situs as this new geometry was called in these days--- to study solutions of differential equations 
by looking at the global structure of the spaces on which they were defined \cite{poincare2010papers}.

This global information was of a different nature than the classical geometrical properties like distances and angles. Unlike the latter, topological
properties of a space do not change under continuous deformations like stretching and bending. They describe more robust features such as the number of holes in a space or 
the different ways to connect to points. Although some of these concepts were already known before, Poincar\'e brought the subject to a new level by demonstrating how they could be used
to derive higly nontrivial results in seemingly unrelated research areas.

This whole new way of thinking formed a fertile playing ground for a new generation of mathematicians. 
Among them was a brilliant Dutch student at the brink of his career and looking for a topic to make his name and fame.
Afther his PhD at the University of Amsterdam, Luitzen Egbert Johannes Brouwer attended the international conference of mathematics in Rome in 1908 \cite{van2012lej}. Impressed 
by Poincar\'e and his ideas, Brouwer started to work on questions on vector fields on surfaces and mappings between spaces, which
resulted in several papers about these topics.

In his papers on vector fields \cite{brouwer1908continuous,freudenthal1976lej} Brouwer showed that every vector field on the sphere has a singular point. 
This is known as the \emph{hairy ball theorem} because in colloquial terms it states that you cannot comb a hairy ball (like a tennis ball) without creating a crown.
This was already proved by Poincar\'e but in a more restricted setting. In his further study Brouwer looked closer at the singular points that can occur by introducing
the winding number of a singularity in a purely topological way.

\def\theX{0.5*(1-cos(deg(x)))} % 0 to 1 to 0
\def\theY{0.5*(sin(deg(x))} % 0 to 0.5 to 0 to -0.5 to 0
\vspace{.3cm}\begin{minipage}{\linewidth}\begin{center}
\begin{tikzpicture}
    % Compute camera unit vector for calculating depth
    \pgfmathsetmacro{\CameraX}{sin(\ViewAzimuth)*cos(\ViewElevation)}
    \pgfmathsetmacro{\CameraY}{-cos(\ViewAzimuth)*cos(\ViewElevation)}
    \pgfmathsetmacro{\CameraZ}{sin(\ViewElevation)}
    \path[use as bounding box] (-1,-1) rectangle (1,1); % Avoid jittering animation
    % Draw a nice looking sphere
    \begin{scope}
        \clip (0,0) circle (1);
        \begin{scope}[transform canvas={rotate=-20}]
            \shade [ball color=white] (0,0.5) ellipse (1.8 and 1.5);
        \end{scope}
    \end{scope}
    \begin{axis}[
        hide axis,
        view={\ViewAzimuth}{\ViewElevation},     % Set view angle
        every axis plot/.style={very thin},
        disabledatascaling,                      % Align PGFPlots coordinates with TikZ
        anchor=origin,                           % Align PGFPlots coordinates with TikZ
        viewport={\ViewAzimuth}{\ViewElevation}, % Align PGFPlots coordinates with TikZ
    ]
\foreach \i in {-5,...,5}{
        \addFGBGplot[domain=0:2*pi, samples=20, samples y=1, black]
(
    {sin(\i*30)*\theX},       % X coordinate
    {2*sin(0.5*\i*30)*\theY}, % Y coordinate
    {1-((1-cos(\i*30)))*\theX} % Z coordinate
);
        }
    \end{axis}
\end{tikzpicture}
\hspace{2cm}
\begin{tikzpicture}
    % Compute camera unit vector for calculating depth
    \pgfmathsetmacro{\CameraX}{sin(\ViewAzimuth)*cos(\ViewElevation)}
    \pgfmathsetmacro{\CameraY}{-cos(\ViewAzimuth)*cos(\ViewElevation)}
    \pgfmathsetmacro{\CameraZ}{sin(\ViewElevation)}
    \path[use as bounding box] (-1,-1) rectangle (1,1); % Avoid jittering animation
    % Draw a nice looking sphere
    \begin{scope}
        \clip (0,0) circle (1);
        \begin{scope}[transform canvas={rotate=-20}]
            \shade [ball color=white] (0,0.5) ellipse (1.8 and 1.5);
        \end{scope}
    \end{scope}
    \begin{axis}[
        hide axis,
        view={\ViewAzimuth}{\ViewElevation},     % Set view angle
        every axis plot/.style={very thin},
        disabledatascaling,                      % Align PGFPlots coordinates with TikZ
        anchor=origin,                           % Align PGFPlots coordinates with TikZ
        viewport={\ViewAzimuth}{\ViewElevation}, % Align PGFPlots coordinates with TikZ
    ]
\foreach \i in {-5,...,5}{
        \addFGBGplot[domain=0:2*pi, samples=20, samples y=1, black]
(
    {cos(deg(x)},       % X coordinate
    {sin(deg(x))*cos(\i*36)}, % Y coordinate
    {sin(deg(x))*sin(\i*36)} % Z coordinate
);
        }
    \end{axis}
\end{tikzpicture}


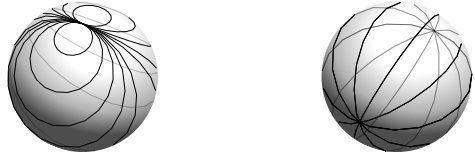
\captionof{figure}{The hairy ball theorem: A vector field on a sphere with one singular point and one with two singular points.}
  \label{fig:hairy}
\vspace{.3cm}\end{center}\end{minipage}

The winding number of a singular point measures how many times the direction of a vector field turns around if one follows a small loop around the singular point.
To make this precise one can look at the map $\omega: \SSS_1  \to \SSS_1$ where the first $\SSS_1$ is the little loop around the singularity, the second
$\SSS_1$ is the set of vector directions in the plane and $\omega$ maps each point in the loop to its vector direction.
If we take a point $p \in \SSS_1$, its preimage $\omega^{-1}(p)$ will be a finite set. For each point $q \in \omega^{-1}(p)$ we can check whether
$\omega$ is orientation preserving at $q$ or not and we can count the difference between the number of orientation preserving points and the orientation reversing points.
For all $p$ for which it is defined, this difference is the same and it is called the winding number.

\vspace{.3cm}\begin{minipage}{\linewidth}\begin{center}
  \begin{tikzpicture}
\draw[-<-] (360*1/7:1)--(0,0); 
\draw[-<-] (432/7:1) .. controls +(1692/7:1) and +(1908/7:1) .. (648/7:1);
\draw[-<-] (72:1) .. controls +(252:.5) and +(1836/7:.5) .. (576/7:1);
\draw[->-] (360*2/7:1)--(0,0); 
\draw[->-] (792/7:1) .. controls +(2052/7:1) and +(324:1) .. (144:1);
\draw[->-] (864/7:1) .. controls +(2124/7:.5) and +(2196/7:.5) .. (936/7:1);
\draw[-<-] (360*3/7:1)--(0,0); 
\draw[-<-] (360*16/35:1)--(0,0); 
\draw[-<-] (360*17/35:1)--(0,0); 
\draw[-<-] (360*18/35:1)--(0,0); 
\draw[-<-] (360*19/35:1)--(0,0); 
\draw[-<-] (360*4/7:1)--(0,0); 
\draw[-<-] (0,0) .. controls +(1728/7:1.3) and +(216:1.3) .. (0,0);
\draw[-<-] (0,0) .. controls +(1656/7:.9) and +(1584/7:.9) .. (0,0);
\draw[->-] (360*5/7:1)--(0,0); 
\draw[->-] (1872/7:1) .. controls +(3132/7:1) and +(3348/7:1) .. (2088/7:1);
\draw[->-] (1944/7:1) .. controls +(3204/7:.5) and +(468:.5) .. (288:1);
\draw[-<-] (360*6/7:1)--(0,0); 
\draw[-<-] (2232/7:1) .. controls +(3492/7:1) and +(3708/7:1) .. (2448/7:1);
\draw[-<-] (2304/7:1) .. controls +(3564/7:.5) and +(3636/7:.5) .. (2376/7:1);
\draw[->-] (360*1:1)--(0,0); 
\draw[->-] (2592/7:1) .. controls +(3852/7:1) and +(4068/7:1) .. (2808/7:1);
\draw[->-] (2664/7:1) .. controls +(3924/7:.5) and +(3996/7:.5) .. (2736/7:1);
\draw[->-,dotted] (30:1.2) arc (30:390:1.2);
\draw (360*1/7:1.2) node[circle,draw,fill=white,minimum size=10pt,inner sep=1pt] {{\tiny1}};
\draw (360*2/7:1.2) node[circle,draw,fill=white,minimum size=10pt,inner sep=1pt] {{\tiny2}};
\draw (360*3/7:1.2) node[circle,draw,fill=white,minimum size=10pt,inner sep=1pt] {{\tiny3}};
\draw (360*4/7:1.2) node[circle,draw,fill=white,minimum size=10pt,inner sep=1pt] {{\tiny4}};
\draw (360*5/7:1.2) node[circle,draw,fill=white,minimum size=10pt,inner sep=1pt] {{\tiny5}};
\draw (360*6/7:1.2) node[circle,draw,fill=white,minimum size=10pt,inner sep=1pt] {{\tiny6}};
\draw (360*7/7:1.2) node[circle,draw,fill=white,minimum size=10pt,inner sep=1pt] {{\tiny0}};

\begin{scope}[xshift=3cm]
\draw[-<-] (360*1/7:1.12)--(360*1/7:.875); 
\draw[-<-] (360*6/35:1.12)--(360*7/35:.875); 
\draw[-<-] (360*7/35:1)--(360*8/35:1); 
\draw[->-] (360*9/35:1.12)--(360*8/35:.875); 
\draw[->-] (360*2/7:1.12)--(360*2/7:.875); 
\draw[->-] (360*11/35:1.12)--(360*12/35:.875); 
\draw[->-] (360*12/35:1)--(360*13/35:1); 
\draw[-<-] (360*14/35:1.12)--(360*13/35:.875); 
\draw[-<-] (360*3/7:1.12)--(360*3/7:.875); 
\draw[-<-] (360*16/35:1.12)--(360*16/35:.875); 
\draw[-<-] (360*17/35:1.12)--(360*17/35:.875); 
\draw[-<-] (360*18/35:1.12)--(360*18/35:.875); 
\draw[-<-] (360*19/35:1.12)--(360*19/35:.875); 
\draw[-<-] (360*4/7:1.12)--(360*4/7:.875); 
\draw[-<-] (360*22/35:1.12)--(360*21/35:.875); 
\draw[->-] (360*22/35:1)--(360*23/35:1); 
\draw[->-] (360*23/35:1.12)--(360*24/35:.875); 
\draw[->-] (360*5/7:1.12)--(360*5/7:.875); 
\draw[->-] (360*26/35:1.12)--(360*27/35:.875); 
\draw[->-] (360*27/35:1)--(360*28/35:1); 
\draw[-<-] (360*29/35:1.12)--(360*28/35:.875); 
\draw[-<-] (360*6/7:1.12)--(360*6/7:.875); 
\draw[-<-] (360*31/35:1.12)--(360*32/35:.875); 
\draw[-<-] (360*32/35:1)--(360*33/35:1); 
\draw[->-] (360*34/35:1.12)--(360*33/35:.875); 
\draw[->-] (360*1:1.12)--(360*7/7:.875); 
\draw[->-] (360*1/35:1.12)--(360*2/35:.875); 
\draw[->-] (360*2/35:1)--(360*3/35:1); 
\draw[-<-] (360*4/35:1.12)--(360*3/35:.875); 
\draw[->-,dotted] (30:1.2) arc (30:390:1.2);
\draw (360*1/7:1.2) node[circle,draw,fill=white,minimum size=10pt,inner sep=1pt] {{\tiny1}};
\draw (360*2/7:1.2) node[circle,draw,fill=white,minimum size=10pt,inner sep=1pt] {{\tiny2}};
\draw (360*3/7:1.2) node[circle,draw,fill=white,minimum size=10pt,inner sep=1pt] {{\tiny3}};
\draw (360*4/7:1.2) node[circle,draw,fill=white,minimum size=10pt,inner sep=1pt] {{\tiny4}};
\draw (360*5/7:1.2) node[circle,draw,fill=white,minimum size=10pt,inner sep=1pt] {{\tiny5}};
\draw (360*6/7:1.2) node[circle,draw,fill=white,minimum size=10pt,inner sep=1pt] {{\tiny6}};
\draw (360*7/7:1.2) node[circle,draw,fill=white,minimum size=10pt,inner sep=1pt] {{\tiny0}};
\end{scope}

\begin{scope}[xshift=6cm]
\draw[->-] (180:1) arc (180:360*1/7:1.1); 
\draw[->-] (360*1/7:1.1) arc (360*1/7:360*2/7-180:1.1);
\draw[->-] (360*2/7-180:1.2) arc (360*2/7-180:360*3/7-360:1.2);
\draw[->-] (360*3/7-360:1.3) arc (360*3/7-360:360*4/7-360:1.3);
\draw[->-] (360*4/7-360:1.4) arc (360*4/7-360:360*5/7-180:1.4);
\draw[->-] (360*5/7-180:1.5) arc (360*5/7-180:360*6/7-360:1.5);
\draw[->-] (360*6/7-360:1.6) arc (360*6/7-360:360*7/7-360-170:1.6)-- (180:1);

\draw (360*1/7:1.1) node[circle,draw,fill=white,minimum size=8pt,inner sep=0pt] {{\tiny1}};
\draw (360*2/7-180:1.2) node[circle,draw,fill=white,minimum size=8pt,inner sep=0pt] {{\tiny2}};
\draw (360*3/7:1.3) node[circle,draw,fill=white,minimum size=8pt,inner sep=0pt] {{\tiny3}};
\draw (360*4/7:1.4) node[circle,draw,fill=white,minimum size=8pt,inner sep=0pt] {{\tiny4}};
\draw (360*5/7-180:1.5) node[circle,draw,fill=white,minimum size=8pt,inner sep=0pt] {{\tiny5}};
\draw (360*6/7:1.6) node[circle,draw,fill=white,minimum size=8pt,inner sep=0pt] {{\tiny6}};
\draw (360*7/7-180:1) node[circle,draw,fill=white,minimum size=8pt,inner sep=0pt] {{\tiny0}};

\draw (-20:.5) arc (-20:20:.5) -- (20:2) arc (20:-20:2) -- (-20:.5);
\draw (0:1.1) node[circle,draw,fill=white,minimum size=8pt,inner sep=0pt] {{\tiny -}};
\draw (0:1.35) node[circle,draw,fill=white,minimum size=8pt,inner sep=0pt] {{\tiny +}};
\draw (0:1.55) node[circle,draw,fill=white,minimum size=8pt,inner sep=0pt] {{\tiny -}};
\draw  (0:2.4) node {{\tiny $=-1$}};

\draw[->-] (-180:.5) arc (-180:180:.5);
\draw[->,thick] (90:.9) -- (90:.6);
\draw (90:.75) node[left] {$\omega$};
\draw (0:.5) node {$\bullet$};
\draw (0:.5) node[left] {$p$};
\end{scope}
\end{tikzpicture}
  
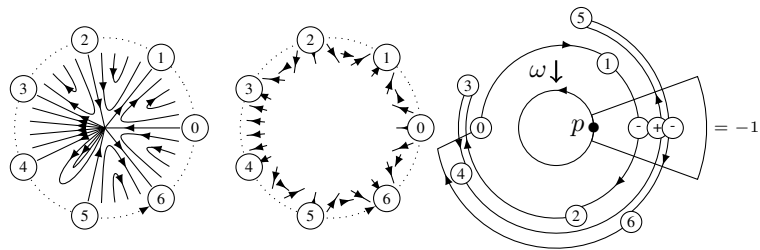
\captionof{figure}{A singularity with winding number $-1$}
  \label{fig:winding}
\vspace{.3cm}\end{center}\end{minipage}

Brouwer noted that this definition made sense for any mapping $\pi:M\to N$ between two oriented manifolds and he called it the mapping degree. 
He realized the power of this new concept and explored its consequences and applications in the paper \emph{\"Uber Abbildung von Mannigfaltigkeiten} \cite{brouwer1911abbildung}. In this paper he not only introduced
the mapping degree but also many other concepts that have become indispensable tools in geometry and topology, such as simplices, simplicial approximations and homotopy of maps.
This paper had an enormous impact on mathematics and it ended with a result ---almost an afterthought--- that would make Brouwer a household name: the Brouwer fixed point theorem.
The fixed point theorem states that every continuous map from the $n$-simplex (or closed $n$-ball) to itself has a point that is fixed by the map. 

\vspace{.3cm}\begin{minipage}{\linewidth}\begin{center}
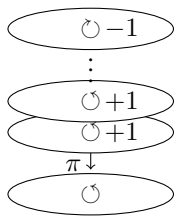

\begin{tikzpicture}[scale=.55]
    \draw[->] (0,1.2)--(0,.1);
    \draw (0,.2) node[left] {$\pi$};
    \draw [fill=white] (0,3.5) ellipse (2 and .5); 
    \draw [fill=white] (0,1) ellipse (2 and .5); 
    \draw [fill=white] (0,1.75) ellipse (2 and .5); 
    \draw [fill=white] (0,-0.5) ellipse (2 and .5); 
    \draw(0,2.7) node {$\vdots$};
    \draw(0,-.5) node {$\circlearrowleft$};
    \draw(0,1) node {$\circlearrowleft$};
    \draw(0,1.75) node {$\circlearrowleft$};
    \draw(0,3.5) node {$\circlearrowright$};
    \draw(.75,1) node {$+1$};
    \draw(.75,1.75) node {$+1$};
    \draw(.75,3.5) node {$-1$};
    \end{tikzpicture}
  \captionof{figure}{The Brouwer degree of a mapping}
  \label{fig:index}
\vspace{.3cm}\end{center}\end{minipage}

The theorem is usually explained in worldly terms by looking at a cup of coffee. In this setting it states that no matter how you stir your cup, there will always be 
a point in the liquid that did not change position and if you try to move that part by further stirring you will inevitably move some other part back into its original position.
Legend even has it that Brouwer came up with the idea while stirring in a real cup, but whether this is true we'll never know. What is true however is that Brouwers refections
on the topic had a profound impact on mathematics and would lead to lots of new developments in geometry. Brouwer himself however moved away from the subject to concentrate
on more foundational issues in mathematics and he left his geometrical ideas for others to explore \cite{van2012lej}. 

\vspace{.3cm}\begin{minipage}{\linewidth}\begin{center}
\resizebox{!}{3cm}{\includegraphics{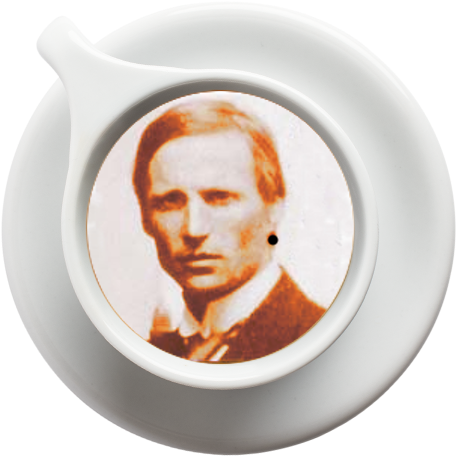}}
\hspace{2cm}
\resizebox{!}{3cm}{\includegraphics{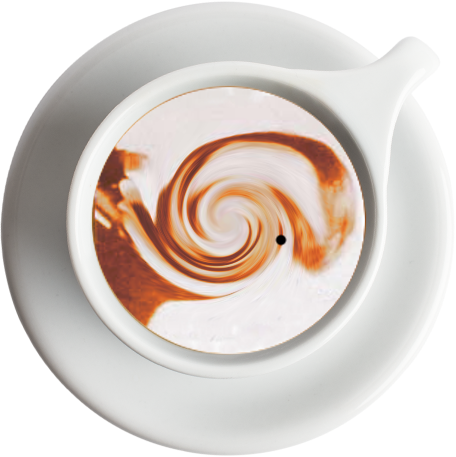}}
  \captionof{figure}{Brouwer's fixed point theorem. The black dot stays place after a rotation and a swirl.}
  \label{fig:swirl}
\vspace{.3cm}\end{center}\end{minipage}

Many new developments in mathematics can be traced back to Brouwer's early investigations and it is an impossible task to discuss them all. Instead we will
highlight a path that takes us from Brouwer's beginnings along some important geometrical theories from the twentieth century, finally leading to mirror symmetry,
a theory that has huge impact on geometry in the early 21st century, just as topology did at the beginning of the 20th century.  

\section{Homology}

The idea that the shape of a space dictates the existence of singularities and fixed points can be refined using a tool that has permeated 20th century mathematics: homology \cite{dieudonne2009history}.
The traces of homology can be found in the work of Euler, who noticed that if one takes the alternating sum of the number of vertices, edges and faces of a polyhedron 
that is sphere shaped one always gets $2$. This holds in general: if one triangulates a space by dividing it up in simplices, the alternating sum of the number of
$i$-dimensional simplices does not depend on the specific triangulation but only on the topology of the space. This number is called the Euler characteristic and
was one of the first topological invariants.

\pgfmathsetmacro{\gr}{(1+sqrt(5))/2}
\pgfmathsetmacro{\igr}{2/(1+sqrt(5))}

%choose axis angles
\newcommand{\xangle}{0}
\newcommand{\yangle}{90}
\newcommand{\zangle}{225}

%choose axis lengths
\newcommand{\xlength}{1}
\newcommand{\ylength}{1}
\newcommand{\zlength}{0.5}

\pgfmathsetmacro{\xx}{\xlength*cos(\xangle)}
\pgfmathsetmacro{\xy}{\xlength*sin(\xangle)}
\pgfmathsetmacro{\yx}{\ylength*cos(\yangle)}
\pgfmathsetmacro{\yy}{\ylength*sin(\yangle)}
\pgfmathsetmacro{\zx}{\zlength*cos(\zangle)}
\pgfmathsetmacro{\zy}{\zlength*sin(\zangle)}

\vspace{.3cm}\begin{minipage}{\linewidth}\begin{center}
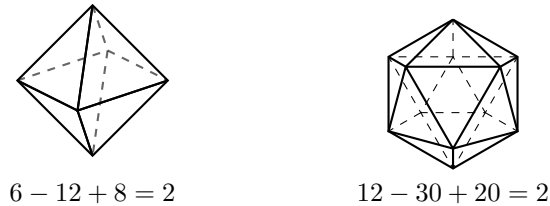

\begin{tikzpicture}[thick,scale=2]
\coordinate (A1) at (0,0);
\coordinate (A2) at (0.6,0.2);
\coordinate (A3) at (1,0);
\coordinate (A4) at (0.4,-0.2);
\coordinate (B1) at (0.5,0.5);
\coordinate (B2) at (0.5,-0.5);

\begin{scope}[thick,dashed,,opacity=0.6]
\draw (A1) -- (A2) -- (A3);
\draw (B1) -- (A2) -- (B2);
\end{scope}
\draw (A1) -- (A4) -- (B1);
\draw (A1) -- (A4) -- (B2);
\draw (A3) -- (A4) -- (B1);
\draw (A3) -- (A4) -- (B2);
\draw (B1) -- (A1) -- (B2) -- (A3) --cycle;
\draw (0.5,-.75) node {$6-12+8=2$};
\end{tikzpicture}
\hspace{2cm}
\begin{tikzpicture}[scale=.33]
    \foreach \y[count=\a] in {10,9,4}
      {\pgfmathtruncatemacro{\kn}{120*\a-90}
       \coordinate (b\a) at (\kn:3);}
    \foreach \y[count=\a] in {8,7,2}
      {\pgfmathtruncatemacro{\kn}{120*\a-90}
       \coordinate (d\a) at (\kn:2.2);}
    \foreach \y[count=\a] in {1,5,6}
      {\pgfmathtruncatemacro{\jn}{120*\a-30}
       \coordinate (a\a) at (\jn:1.5);}
    \foreach \y[count=\a] in {3,11,12}
      {\pgfmathtruncatemacro{\jn}{120*\a-30}
       \coordinate (c\a) at (\jn:3);}
  \draw[dashed] (a1)--(a2)--(a3)--(a1);
  \draw[thick] (d1)--(d2)--(d3)--(d1);
  \foreach \a in {1,2,3}
   {\draw[dashed] (a\a)--(c\a);
   \draw[thick] (d\a)--(b\a);}
   \draw[thick] (c1)--(b1)--(c3)--(b3)--(c2)--(b2)--(c1);
   \draw[thick] (c1)--(d1)--(c3)--(d3)--(c2)--(d2)--(c1);
   \draw[dashed] (b1)--(a1)--(b2)--(a2)--(b3)--(a3)--(b1);
\draw (0,-4) node {$12-30+20=2$};
 \end{tikzpicture}
  \captionof{figure}{The Euler characteristic}
  \label{fig:Euler}

  \vspace{.3cm}\end{center}\end{minipage}

The Euler characteristic is a powerful invariant that pops up in a lot of situations. In our story Poincar\'e used it to refine the existence theorem of singularities of vector 
fields on surfaces. He showed that the sum of all the winding numbers of the singular points equals the Euler characteristic of the surface. Because the Euler characteristic
of the sphere is $2$, there must be at least one singular point for a vector field on a sphere.

This result was generalized by Hopf to higher dimensions. Using Brouwer's idea of the mapping degree one can define the index of a singular point of a vector field in higher
dimensions. Just like the winding number in two dimensions, it is the degree of the map that maps a small $n$-sphere around the singularity to the $n$-sphere of vector 
directions. Hopf proved that the sum of the indices of the singular points of a vector field always equals the Euler characteristic of the space. 
This result is known as the Poincar\'e-Hopf theorem \cite{guillemin2010differential}.

To get more detailed results about the singularities of a vector field, one needs a refinement of the Euler characteristic. 
The Betti numbers, introduced by Betti and put on a firmer grounding by Poincar\'e \cite{poincare2010papers}, measure the number of holes in a space. 
Their idea was that holes can be detected by looking at boundaries. If you look at the boundary of a submanifold, you may notice that it has no boundary itself: 
e.g. the boundary of a disk is a circle and a circle has no end points. A subspace without boundary is called a cycle, so every boundary is a cycle.
The converse is however not always true, some cycles are not the boundary of 
something else e.g. if we look at a circle embedded in a ring this circle is not the boundary of something because there is a hole in the ring.
Poincar\'e defined the $n^{th}$-Betti number as the maximal number of $n$-dimensional cycles that do not form the boundary of an $n+1$-dimensional submanifold.

To turn this into something easy to work with, it is best to linearize this construction by introducing a homology theory. The idea is to associate 
to the space a complex. This is a graded vector space $V_\bullet$ with a linear operator $d:V_{\bullet}\to V_{\bullet}$ of degree $-1$, which is called the boundary operater.
Because the boundary of a boundary is empty, we demand that $d^2=0$. In a complex we have that $\Image d \subset \Ker d$ and we can calculate the quotient 
$$
H_\bullet(V):= \frac{\Ker d}{\Image d}.
$$
This is called the homology of the complex.

Starting from a triangulation of a manifold $\MM$, we take $V_i$ to be the vector space spanned by all $i$-dimensional simplices and $d$ maps
an $i$-simplex to the signed sum of the $i-1$-simplices that make up its boundary. Although $V_\bullet$ depends highly on the triangulation of the space,
the homology does not: different triangulations give homologies that are isomorphic as graded vector spaces. The dimensions of the $n$-th degree part
of the homology is called the $n$-th Betti number of the space. One can show that the alternating sum of the Betti numbers equals the Euler characteristic and
hence they can be seen as a refinement of this invariant.

\vspace{.3cm}\begin{minipage}{\linewidth}\begin{center}
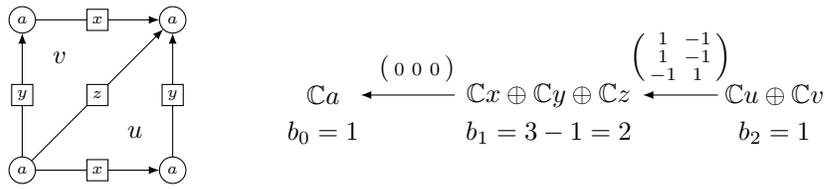


\begin{tikzpicture}[scale=1]
\draw[-latex,shorten >=5pt] (0,0)--(2,0);
\draw[-latex,shorten >=5pt] (2,0)--(2,2);
\draw[-latex,shorten >=5pt] (0,0)--(0,2);
\draw[-latex,shorten >=5pt] (0,2)--(2,2);
\draw[-latex,shorten >=5pt] (0,0)--(2,2);
\draw (0,0) node[circle,draw,fill=white,minimum size=10pt,inner sep=1pt] {{\tiny $a$}};
\draw (2,0) node[circle,draw,fill=white,minimum size=10pt,inner sep=1pt] {{\tiny $a$}};
\draw (0,2) node[circle,draw,fill=white,minimum size=10pt,inner sep=1pt] {{\tiny $a$}};
\draw (2,2) node[circle,draw,fill=white,minimum size=10pt,inner sep=1pt] {{\tiny $a$}};
\draw (1,0) node[draw,fill=white,minimum size=8pt,inner sep=1pt] {{\tiny $x$}};
\draw (1,2) node[draw,fill=white,minimum size=8pt,inner sep=1pt] {{\tiny $x$}};
\draw (0,1) node[draw,fill=white,minimum size=8pt,inner sep=1pt] {{\tiny $y$}};
\draw (2,1) node[draw,fill=white,minimum size=8pt,inner sep=1pt] {{\tiny $y$}};
\draw (1,1) node[draw,fill=white,minimum size=8pt,inner sep=1pt] {{\tiny $z$}};
\draw (1.5,.5) node{$u$};
\draw (.5,1.5) node{$v$};

\begin{scope}[xshift=4cm,yshift=1cm]
\draw (0,0) node {$\C a$};
\draw[latex-] (.5,0) -- (1.75,0);
\draw (3,0) node {$\C x \oplus \C y \oplus \C z$};
\draw[latex-] (4.25,0) -- (5.25,0);
\draw (6,0) node {$\C u \oplus \C v$};
\draw (1.25,.35) node{$\sm{0&0&0}$};
\draw (4.75,.5) node{$\sm{1&-1\\1&-1\\-1&1}$};
\draw (0,-.5) node {$b_0=1$};
\draw (3,-.5) node {$b_1=3-1=2$};
\draw (6,-.5) node {$b_2=1$};
\end{scope}
\end{tikzpicture}
  \captionof{figure}{The Betti numbers of a torus.}
  \label{fig:betti}
\vspace{.3cm}\end{center}\end{minipage}

A triangulation is just one way to cook up a complex. This technique goes under the name of simplicial homology, but
soon there were many others using different ideas in geometry: cellular homology, singular homology, de Rham cohomology, etc.
All gave different complexes with the same underlying homology and therefore they could be used to explore connections between many
different aspects of geometry \cite{dieudonne2009history}.
This proved to be so useful that it soon spread over the many other fields in mathematics and these days almost all areas in mathematics make use of homology theories.
Homology has become a versatile tool that can be used to measure holes, intersections, deformations and obstructions for many different mathematical objects \cite{weibel1999history}.

We will now split our story in two and highlight two examples of homology theory that arose out of this in the second half of the twentieth century: sheaf cohomology and Floer homology.
We will discuss both theories from the viewpoint of Brouwer's work and then we will braid them together in the final act of our story. 

\section{Sheaves and their cohomology}

Jean Leray was a French mathematician who had worked on extending Brouwer's fixed point theorem to infinite dimensional settings and had used it to prove
existence results for partial differential equations. However by 1940 he found himself locked up as a prisoner of war in a camp in Austria and to bide his time he decided to
devote himself to the study of algebraic topology. During the next five years he developed three ideas that would be very influential after the war: sheaves, sheaf cohomology
and spectral sequences \cite{miller2000leray}.

Geometrically a sheaf over a space $\MM$ can be seen as a generalization of a bundle, from the point of view of its sections. A section of a bundle $\pi : B\to \MM$
assigns in a continuous way to each point $m \in \MM$ an element of its fiber $\pi^{-1}(m)$. If this is done for the whole of $\MM$ then we call this a global section, but if 
it only works for an open part of $\MM$ we call it a local section. Local sections that overlap can be glued together to make bigger local sections and if we formalize
this behaviour we arrive at the notion of a sheaf. 

Brouwer's theorem on vector fields can be restated in sheaf terminology: the sheaf corresponding to the tangent bundle $\T\SSS_2$ does not admit a nowhere vanishing global section.
The existence of global sections and the obstructions in the gluing process are governed by a homology theory, sheaf cohomology, which assigns to each sheaf 
a sequence of homology groups. The classical homology groups of a manifold can be recovered by looking
at the sheaf cohomology of a special sheaf: the constant sheaf with coefficients in $\R$.

\vspace{.3cm}\begin{minipage}{\linewidth}\begin{center}
\begin{tikzpicture}[yscale=.66]
\begin{scope}
    \draw  (0,0) ellipse (2 and .5); 
    \draw [dashed] (0,0.5) ellipse (2 and .5); 
    \draw [dotted] (0,-1.25) ellipse (2 and .5); 
    \draw [dotted] (0,1.25) ellipse (2 and .5); 
    \draw [dotted](2,-1.25)--(2,1.25);
    \draw [dotted](-2,-1.25)--(-2,1.25);
\end{scope}
\begin{scope}[xshift=6cm]
    \draw [dotted](2,-1.25)--(2,1.25);
    \draw [dotted](-2,-1.25)--(-2,1.25);
    \draw (0,0) ellipse (2 and .5); 
    \draw [dotted](2,-1.25) .. controls (2,-.75) and (.5,-.75) .. (0,-.75); 
    \draw [dotted](-2,-1.25) .. controls (-2,-.75) and (-.5,-.75) .. (0,-.75); 
    \draw [dotted](2,-1.25) .. controls (2,-1.75) and (1.5,-1.75) .. (1,-1.75); 
    \draw [dotted](-2,-1.25) .. controls (-2,-1.75) and (-1.5,-1.75) .. (-1,-1.75); 
  \begin{scope}[yshift=2.5cm]    
    \draw [dotted](2,-1.25) .. controls (2,-.75) and (.5,-.75) .. (0,-.75); 
    \draw [dotted](-2,-1.25) .. controls (-2,-.75) and (-.5,-.75) .. (0,-.75); 
    \draw [dotted](2,-1.25) .. controls (2,-1.75) and (1.5,-1.75) .. (1,-1.75); 
    \draw [dotted](-2,-1.25) .. controls (-2,-1.75) and (-1.5,-1.75) .. (-1,-1.75);  
  \end{scope}
  \draw[dotted] (-1,.75) .. controls (-.5,.75) and (.5,-1.75) ..(1,-1.75);
  \draw[dotted] (-1,.75) .. controls (-.5,.75) and (.5,-1.75) ..(1,-1.75);
  \draw[dotted] (-1,-1.75) .. controls (-.5,-1.75) and (-.25,-1.05) ..(-0.05,.-.6);
    \draw[dotted] (0.05,-.40) .. controls (.25,.05) and (.5,.75) ..(1,.75);
\begin{scope}[yshift=-.25cm]
\begin{scope}[yscale=.55]
  \draw[dashed] (-1,.75) .. controls (-.5,.75) and (.5,-1.75) ..(1,-1.75);
  \draw[dashed] (-1,.75) .. controls (-.5,.75) and (.5,-1.75) ..(1,-1.75);
  \draw[dashed] (-1,-1.75) .. controls (-.5,-1.75) and (-.25,-1.05) ..(-0.05,.-.6);
    \draw[dashed] (0.05,-.40) .. controls (.25,.05) and (.5,.75) ..(1,.75);
    \draw[dashed] (2,-1.25) .. controls (2,-.75) and (.5,-.75) .. (0,-.75); 
    \draw[dashed] (-2,-1.25) .. controls (-2,-.75) and (-.5,-.75) .. (0,-.75); 
    \draw[dashed] (2,-1.25) .. controls (2,-1.75) and (1.5,-1.75) .. (1,-1.75); 
    \draw[dashed] (-2,-1.25) .. controls (-2,-1.75) and (-1.5,-1.75) .. (-1,-1.75); 
  \begin{scope}[yshift=2.5cm]    
    \draw[dashed] (2,-1.25) .. controls (2,-.75) and (.5,-.75) .. (0,-.75); 
    \draw[dashed] (-2,-1.25) .. controls (-2,-.75) and (-.5,-.75) .. (0,-.75); 
    \draw[dashed] (2,-1.25) .. controls (2,-1.75) and (1.5,-1.75) .. (1,-1.75); 
    \draw[dashed] (-2,-1.25) .. controls (-2,-1.75) and (-1.5,-1.75) .. (-1,-1.75);  
  \end{scope}
  \draw[dashed] (-1,.75) .. controls (-.5,.75) and (.5,-1.75) ..(1,-1.75);
  \draw[dashed] (-1,.75) .. controls (-.5,.75) and (.5,-1.75) ..(1,-1.75);
  \draw[dashed] (-1,-1.75) .. controls (-.5,-1.75) and (-.25,-1.05) ..(-0.05,.-.6);
    \draw[dashed] (0.05,-.40) .. controls (.25,.05) and (.5,.75) ..(1,.75);
\end{scope}
\end{scope}
    %   \draw (0,1) ellipse (2 and .5); 
\end{scope}

%    .. controls (1,0) and (0,-1) .. (-2,-1) .. controls (0,1) and (1,0) .. (1,-1.5); 
%    \draw[yshift=2] (3,-1) .. controls (0,1) and (-1,0) .. (0,-.5) .. controls (1,0) and (0,-1) .. (-2,-1) .. controls (0,1) and (1,0) .. (1,-1.5); 
\end{tikzpicture}
  
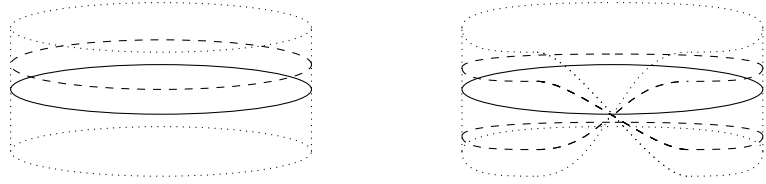
\captionof{figure}{The cylinder sheaf over the circle admits a nonvanishing section while the Moebius sheaf does not}
  \label{fig:sheaf}
\vspace{.3cm}\end{center}\end{minipage}

After returning to France, Leray published his wartime exploits in the Comptes Rendues where they were picked up by Cartan and Serre and many other
mathematicians in Paris. Inspired by these homological techniques, Paris in the fifties and sixties became the place and time for a revolution in complex and algebraic
geometry. Under the lead of Alexander Grothendieck the field was reworked, reinterpreted and extended beyond recognition. The new algebraic geometry was firmly
founded on the notions of sheaves and their different cohomology theories. These ideas were used to prove many big results such as a far-reaching generalization of the Riemann-Roch theorem and 
the celebrated Weil conjectures \cite{dieudonne1985history}.

We can illustrate the idea of sheaf cohomology with an example. 
The Riemann sphere $\PP_1$ is an 
algebraic variety that can be constructed by gluing two copies of $\C$ together by identifying the coordinate
of the first copy with the inverse of the coordinate of the second copy. The polynomial
functions on the two parts can be seen as $\C[z]$ and $\C[z^{-1}]$. A sheaf $\sM$ on $\PP_1$ consists
of a $\C[z]$-module $M_1$, a $\C[z^{-1}]$-module $M_2$ and a $\C[z,z^{-1}]$-module $M_{12}$ 
together with restriction morphisms $\phi_i:M_i \to M_{12}$.
A global section of $\sM$ consists of a pair $(m_1,m_2) \in M_1\oplus M_2$ that agrees on the intersection:
$\phi_1(m_1)=\phi_2(m_2)$. These can be seen as elements in the zeroth cohomology of the following complex.
\[
\begin{tikzpicture}
 \draw (-0.5,0) node{$M_1\oplus M_2$};
 \draw (1.5,0.3) node{$\phi_1 \oplus -\phi_2$};
 \draw (3,0) node{$M_{12}$};
 \draw[-latex] (.5,0)--(2.5,0);
\end{tikzpicture}
\]
The elements of the first cohomology are elements in $M_{12}$ that do not come
from restrictions of elements of $M_1$ or $M_2$. In other words these are local sections
that are obstructed. 

An easy example of a sheaf is $\sO(i)$. Here $M_{12} = \C[z,z^{-1}]$, $M_1=\C[z]$, $M_2=z^i\C[z^{-1}]$ and 
the restriction maps are embeddings. From the discussion above it
is clear that $\sO(i)$ has $1+i$ linearly independent global sections if $i\ge 0$ and 
none otherwise. On the other hand the first cohomology group of $\sO(i)$ is $-1-i$-dimensional
if $i< 0$ and zero otherwise:
\[
 \dim H^0(\sO(i))=\begin{cases}
                   i+1&i\ge 0\\
                   0&i< 0
                  \end{cases}
\hspace{1cm}
 \dim H^1(\sO(i))=\begin{cases}
                   0&i\ge 0\\
                   -1-i&i< 0.
                  \end{cases}
                \]

\section{Morse and Floer}

Morse theory was first developed by Marston Morse \cite{bott1988morse}, an American professor with whom Brouwer corresponded and whom he also payed a visit at the Institute
of Advanced Studies in Princeton \cite{van2012lej}. Morse's original idea was to study the shape of a manifold by slicing it like a boiled egg and studying how 
the slices differ from another.
The slicing tool he used was a Morse function $f:M \to \R$, which is a smooth function with nondegenerate singularities. The singularities of $f$ are points where $\left(\frac{ f}{\partial x_i}\right)$ vanishes
and they are nondegenerate if the Hessian $\left(\frac{\partial^2f}{\partial x_i\partial x_j}\right)$ is an invertible matrix. Singularities come in different flavours 
depending on the number
of negative eigenvalues of the Hessian. This number is called the Morse index. Whenever the slicer crosses a singular point, the topology of the slice will change. Morse studied these changes
to find handle decompositions of manifolds.

\vspace{.3cm}\begin{minipage}{\linewidth}\begin{center}
\resizebox{!}{3cm}{\includegraphics{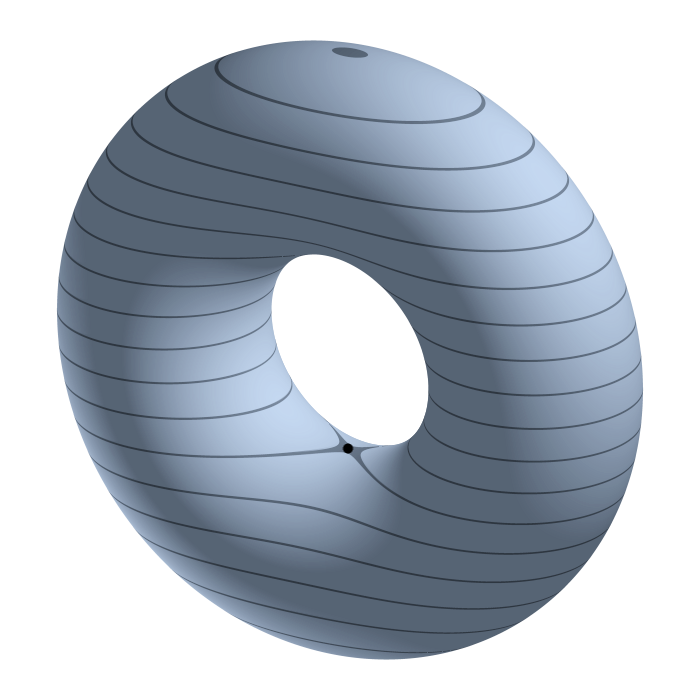}}
\captionof{figure}{Morse theory on a torus: At the black dot, a singular point with Morse index 1, the slice splits in two.}
  \label{fig:slice}
\vspace{.3cm}\end{center}\end{minipage}

The original work of Morse did not emphasize the relation with homology groups. His ideas were revisited by many mathematicians to make this connection more explicit, such as
Milnor, Thom and Smale \cite{hutchings2002lecture}. Finally in 1982 Witten \cite{witten1982supersymmetry} gave this new approach more visibility to the outside world via his work on supersymmetry. The idea is to use a metric 
on the manifold to turn the Morse function $f$ into a vector field by taking the gradient $\nabla f$. Singular points in the Morse function will correspond to singularities of the vector field 
and the Morse index will count the number of incoming directions of the vector field. It is important to note that the Morse index is not the same as the Brouwer index 
for singularities of vector fields. For singularities coming from Morse functions the Brouwer index is $-1$ to the power of the Morse index.

From the vector field $\nabla f$ one can construct a complex. The $i$-degree component is the vector space spanned by the singularities with Morse index $i$ and the differential
$d$ maps each singularity with Morse index $i$ to a sum of all the index $i-1$ singularities that can be reached by an integral curve of the vector field. 
Multiple connecting paths contribute each with a sign. One of the main results in this revisited version of Morse theory is that the homology of this complex equals the simplicial
homology of the space. This result tells us that the number of critical points of $f$ must be at least equal to the sum of all Betti numbers, otherwise the complex could not 
have the required homology. Likewise we have a lower bound on the number of singularities in the vector field $\nabla f$ and the number of fixpoints of the map $\exp(t\nabla f): \MM\to \MM$.

\vspace{.3cm}\begin{minipage}{\linewidth}\begin{center}
\begin{tikzpicture}
\draw (0,0) node {\resizebox{!}{3cm}{\includegraphics{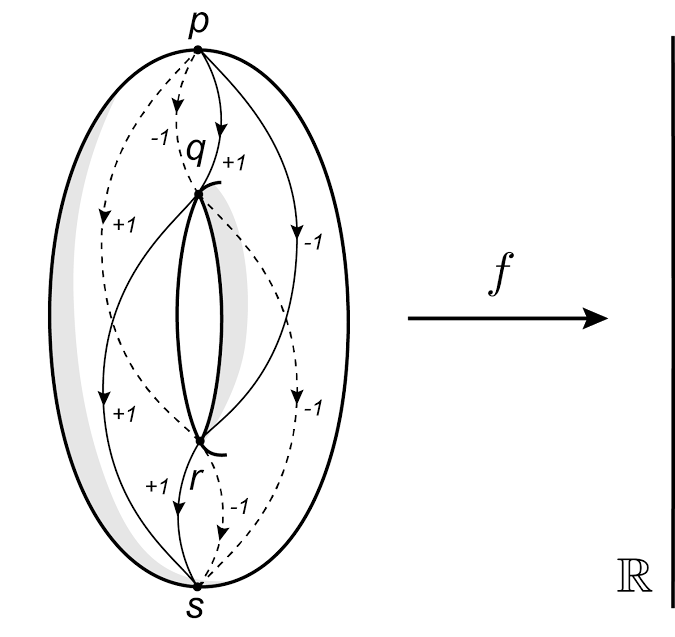}}};
\begin{scope}[xscale=.75, xshift=4cm,yshift=.5cm]
\draw (0,0) node {$\C s$};
\draw[latex-] (.5,0) -- (1.75,0);
\draw (3,0) node {$\C q \oplus \C r$};
\draw[latex-] (4.25,0) -- (5.25,0);
\draw (6,0) node {$\C p$};
\draw (1.25,.35) node{$\sm{0&0}$};
\draw (4.75,.5) node{$\sm{0\\0}$};
\draw (0,-.5) node {$b_0=1$};
\draw (3,-.5) node {$b_1=2$};
\draw (6,-.5) node {$b_2=1$};
\end{scope}
\end{tikzpicture}
  \captionof{figure}{The Morse complex of a torus: each path is canceled by a path with opposite sign, so all differentials are zero.}
  (credit: Jeremy van der Heyden)
  \label{fig:comp}
\vspace{.3cm}\end{center}\end{minipage}

The Morse complex can also be given a more physical interpretation: the manifold $\MM$ represents the space in which a particle can move under a force given by a potential $f$. 
The singularities are the points where the particle is in equilibrium and the Morse index measures how stable these equilibria are. The differential counts how many
paths there are between the equilibria.
This dynamical interpretation of Morse theory indicated that this theory could somehow be extended to more general dynamic settings. In mathematical terms this meant
an excursion into symplectic geometry. 

Symplectic manifolds were introduced by Hermann Weyl as a mathematical model for phase space in physics. They are even-dimensional manifolds equiped with a closed two-form, 
that can be used to measure areas. The two-form also allows us to identify the tangent bundle and the cotangent bundle, so the differential $df$ of a function
$f$ can be turned into a vector field called the Hamiltonian vector field. This vector field gives a dynamical flow: the Hamiltonian flow. 
The standard example is the cotangent space of a manifold $\T^*\MM$ with the two-form $\omega = \sum dp_i\wedge dq_i$ and if we take for $f$ the total energy of the system, 
the Hamiltonian flow will give the dynamics of the physical system.

If we go back to the example from Morse homology, the function $f: M \to \R$ can be turned into a Hamiltonian
by adding a kinetic energy term. The Hamiltonian vector field will project onto $\nabla f$ and hence describes our physical model of Morse theory. 
The fixed points of the Hamiltonian vector field are the singular points of $f$, but they can also be interpreted as intersection points.
The manifold $\MM$ sits inside $\T^*\MM$ as a maximal-dimensional subspace on which $\omega$ vanishes. Such subspaces are called Lagrangian subspaces and a Hamiltonian flow
maps a Lagrangian subspace to a Lagrangian subspace. If we look at the flowed version of the $\phi_t(\MM)\subset \T^*\MM$ we see that it intersects $\MM$ at the points
where $f$ is singular, so the Morse complex of $(M,f)$ is spanned by the intersection points of  $\phi_t(\MM)$ and $\MM$ and the connecting paths swipe out strips
with boundaries on the two submanifolds

\vspace{.3cm}\begin{minipage}{\linewidth}\begin{center}
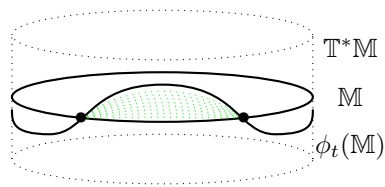

\begin{tikzpicture}[yscale=.66]
\begin{scope}
    \draw [thick] (0,0) ellipse (2 and .5); 
    \draw (2.5,0) node {$\MM$};
    \draw (2.5,1) node {$\T^*\MM$};
    \draw (2.5,-1) node {$\phi_t(\MM)$};
    \draw[dotted,green] (-1.08,-.43).. controls (-.5,.3) and (.5,.3) .. (1.08,-.43);
    \draw[dotted,green] (-1.08,-.43).. controls (-.5,.2) and (.5,.2) .. (1.08,-.43);
    \draw[dotted,green] (-1.08,-.43).. controls (-.5,.1) and (.5,.1) .. (1.08,-.43);
    \draw[dotted,green] (-1.08,-.43).. controls (-.5,0) and (.5,0) .. (1.08,-.43);
    \draw[dotted,green] (-1.08,-.43).. controls (-.5,-.1) and (.5,-.1) .. (1.08,-.43);
    \draw[dotted,green] (-1.08,-.43).. controls (-.5,-.2) and (.5,-.2) .. (1.08,-.43);
    \draw[dotted,green] (-1.08,-.43).. controls (-.5,-.3) and (.5,-.3) .. (1.08,-.43);
    \draw[dotted,green] (-1.08,-.43).. controls (-.5,-.4) and (.5,-.4) .. (1.08,-.43);
    \draw [thick] (2,-.25) .. controls (2,-.5) and (2,-.75) .. (1.5,-.75) .. controls (1,-.75) and (1,.25) .. (0,.25) .. controls (-1,.25) and (-1,-.75) .. (-1.5,-.75) .. controls (-2,-.75) and (-2,-.5) .. (-2,-.25); ; 
    %\draw [dotted] (0,0.5) ellipse (2 and .5); 
    \draw [dotted](0,-1.25) ellipse (2 and .5); 
    \draw [dotted] (0,1.25) ellipse (2 and .5); 
    \draw [dotted](2,-1.25)--(2,1.25);
    \draw [dotted] (-2,-1.25)--(-2,1.25);
    \draw (-1.08,-.43) node{$\bullet$};
    \draw (1.08,-.43) node{$\bullet$};
\end{scope}
\end{tikzpicture}
  \captionof{figure}{Morse homology as Floer homology: each critical point induces an intersection point and each path in $\MM$ induces a strip in $\T^*\MM$.}
  \label{fig:floer}
\vspace{.3cm}\end{center}\end{minipage}

This observation lead Andreas Floer to construct a complex for any pair of Lagrangian submanifolds that intersect normally \cite{floer1988morse}. Again the complex is
spanned by the intersection points and the differential is given by strips between them. 
In the Morse case the connecting lines depended on the metric and the Morse function, but ultimately the homology only depended topology of $\MM$.
Similarly Floer needed a (almost)-complex structure on $\MM$ to define his strips as (pseudo)-holomorphic disks, but ultimately the homology of his complex 
only depended on the symplectic topology of the space and the Lagrangian submanifolds. Moreover the homology did not change under Hamiltion flows of 
the Lagrangian submanifolds, so it only depended on the Hamiltonian isotopy classes of the submanifolds.

In the spirit of Brouwer, Floer used his ideas to tackle a famous question about fixed points: Arnold's conjecture \cite{floer1986proof,liu1998floer}. This conjecture states that any
Hamiltonian map of a manifold $\MM$ onto itself has at least as many fixed points as the sum of the Betti numbers of $\MM$. Floer was able to prove this conjecture under certain
conditions. The techniques introduced by Floer opened up a whole new branch in mathematics: symplectic topology. This research area gained a lot of interest in the next decades and 
is still very active today. Sadly Floer did not live to see his legacy, he commited suicide in 1991 \cite{hofer2012floer}.
Coincidentally this was the same year that the final ingredient of our story entered the mathematical scene: mirror symmetry.

\section{The new beginning}

Mirror symmetry has its origins in superstring theory, a theory in quantum physics that aims to unite all forces of nature by postulating that elementary particles are 
not points but tiny strings that move in a ten dimensional world. The ten dimensions consist of the four classical spacetime dimensions and six that are rolled up to form
a compact manifold. This manifold $\MM$ has a lot of extra structure: it has a symplectic structure, a complex structure, and it is Calabi-Yau, which is the complex analogon of being orientable.
The shape of the space determines how the particles behave in spacetime but the complexity of the theory keeps a full study of the implications of the theory currently 
out of reach. 

In order to make sense of the theory, physicists look at approximations that concentrate only certain aspects and this resulted in two models that were being pursued.
The A-model looks at the symplectic structure but ignores the complex structure and the B-model does the opposite \cite{vonk2005mini}. While doing calculations, physicists noticed
that calculations for the A-model of a certain manifold $\MM$ yielded the same results as calculations for the B-model of a different manifold $\MM'$. These two manifolds were called mirrors of each other
because their Hodge numbers, a refinement of the Betti numbers, were reversed. This observation was very useful because in some cases calculations in the A-model were easier,
while in other cases calculations in the B-model were easier. 

At that time, the mid eighties, mathematicians took no notice of this but this changed when in 1991 four physicists, Candelas, De la Ossa, Green and Parkes, \cite{candelas1991pair} used these ideas to conjecture
a formula that counted the number of curves of any degree in the quintic threefold. At that time mathematicians had only counted these curves in very low degree and had 
absolutely no idea what a general formula would look like. After some period of disbelief mathematicians became convinced that this formula really worked but they could not
understand the physical reasoning behind it. Many mathematicians embarked on a quest to uncover the mathematics that made this work. By the turn of the century several mathematical proofs of this formula were found \cite{givental1996equivariant,lian1998mirror}.

In 1994 Kontsevich gave a visionary talk at the International Congress of Mathematicians in Z\"urich \cite{kontsevich1994homological}. In this talk he layed out 
a research program that is now known under the 
name of Homological Mirror Symmetry. He proposed that mirror symmetry should be seen as a correspondence between the homology theories associated to the two
types of geometry involved: Floer theory for symplectic geometry and sheaf cohomology for complex geometry. For each pair of Lagrangian submanifolds in $\MM$ there should exist
a pair of sheaves on the mirror $\MM'$ such that the Floer homology of the Lagrangians is isomorphic to the sheaf homology between the two sheaves\footnote{Technically we need to take complexes of sheaves and Lagrangians}. 
To make this precise he made use of the language of categories. To each of the manifolds one can assign a category that packages the homological informations in a neat way: the derived Fukaya category $\mathtt{DFuk} \MM$ on the symplectic side
and the derived category of coherent sheaves $\mathtt{DCoh} \MM'$. Mirror symmetry manifests itself as an equivalence between these two categories.

The first instance of homological mirror symmetry was worked out by Polishchuk and Zaslow in 1998 \cite{polishchuk2001categorical}. They consider the case where $\MM$ is a symplectic torus and in that case $\MM'$ is also a torus (but with the structure of an elliptic curve).
The Lagrangian submanifolds are closed curves on the torus and they can be characterized by two numbers $(a,b)$ which denote the direction of the line. These numbers must be integers
to make sure that the line closes in on itself on the torus. Sheaves on an elliptic curve have also two integers associated to them: their rank and their degree. By 
identifying these (and other) numerical invariants Polishchuk and Zaslow were indeed able to construct an equivalence of categories.

\vspace{.3cm}\begin{minipage}{\linewidth}\begin{center}
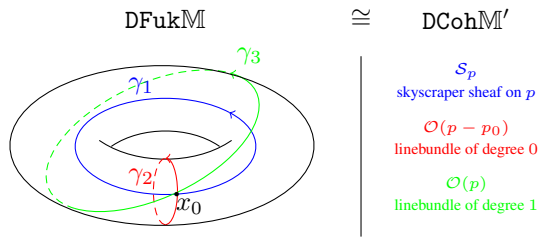

\begin{tikzpicture}  
%Hole
\draw (0,1.7) node {$\mathtt{DFuk} \MM$};
\draw (2.6,1.2)--(2.6,-1.2);
\draw (2.6,1.7) node {$\cong$};
\begin{scope}[scale=.8]
%\begin{scope}[scale=0.8]
\path[rounded corners=24pt] (-.9,0)--(0,.6)--(.9,0) (-.9,0)--(0,-.56)--    (.9,0);
\draw[rounded corners=28pt] (-1.1,.1)--(0,-.6)--(1.1,.1);
\draw[rounded corners=24pt] (-.9,0)--(0,.6)--(.9,0);
%\end{scope}
%Cut 1
\draw[densely dashed,red] (0,-1.3) arc (270:90:.2 and 0.550);
\draw[red,->] (0,-1.3) arc (-90:90:.2 and .550);
%Cut 2
\draw[blue,->] (1.0605,0.5656) arc (45:405:1.5 and 0.8);
\draw[green,rotate=25,->] (-1.8,-.3) arc (210:390:1.9 and .95);
\draw[green,rotate=25,densely dashed] (-1.8,-.3) arc (210:30:1.9 and .95);
\draw[black,scale=1.63] (1.0605,0.5856) arc (45:405:1.55 and 0.815);
\node at (0.19,-0.79) {$\textbf{.}$};
\node at (0.4,-1) {$x_0$};
\node at (-0.4,1) {$\textcolor{blue}{\gamma_1}$};
\node at (-0.4,-0.5) {$\textcolor{red}{\gamma_2}$};
\node at (1.4,1.5) {$\textcolor{green}{\gamma_3}$};
\end{scope}
\draw (4,1.7) node {$\mathtt{DCoh} \MM'$};
\draw (4,1) node[blue] {\tiny $\cS_p$};
\draw (4,.7) node[blue] {\tiny skyscraper sheaf on $p$};
\draw (4,.25) node[red] {\tiny $\cO(p-p_0)$};
\draw (4,-.05) node[red] {\tiny linebundle of degree $0$};
\draw (4,-.5) node[green] {\tiny $\cO(p)$};
\draw (4,-.8) node[green] {\tiny linebundle of degree $1$};

\end{tikzpicture}
  \captionof{figure}{Mirror symmetry for the torus: closed curves on the torus correspond to sheaves on the elliptic curve.}
  \label{fig:mirror}
\vspace{.3cm}\end{center}\end{minipage}

In the new millenium Homological Mirror Symmetry became the inspiration for a lot of mathematical research. Proofs were given in other special instances and the conjecture was 
extended to situations where the mirror is not a manifold but something more exotic such as a singular or noncommutative space. Techniques that were developed to overcome 
difficulties in the study of mirror symmetry had a major impact on the research in related areas in geometry, algebra, category theory and theoretical physics.
This has lead to a new way of thinking about the relations between these different subjects that will influence future of mathematics for many years to come.

\section{Further Reading}

I sketched a path that connects Brouwer's topological work with modern developments in geometry. The path is only one of many interesting walks one can take along the history 
of 20th century geometry and topology. For the interested reader I provide some interesting travel guides along the way. First of all, for those interested in Brouwer's life and legacy I recommend
the biography \emph{L. E. J. Brouwer: Topologist, Intuitionist, Philosopher. How Mathematics is Rooted in Life} by Van Dalen \cite{van2012lej}. More about the history of
algebraic topology and algebraic geometry can be found in the books of Dieudonne, \emph{A History of Algebraic and Differential Topology, 1900 - 1960}, \cite{dieudonne2009history}
 and \emph{History of Algebraic geometry} \cite{dieudonne1985history}.  
A good overview of Morse theory can be found in Bott's article \emph{Morse indomitable} \cite{bott1988morse} and for 
Floer theory I refer to the notes by Salamon \cite{salamon1999lectures}. How Floer theory gave rise to the Fukaya category is treated in Smith's \emph{A symplectic prolegomenon}
\cite{smith2015symplectic}. For the B-model I recommend Cald$\breve{\mathrm{a}}$r$\breve{\mathrm{a}}$ru's lecture notes \emph{Derived categories of sheaves: a skimming} \cite{cualduararu2005derived}.
For homological mirror symmetry I do not know a good historical overview article but several basic examples are worked out in
Ballard's \emph{Meet homological mirror symmetry} \cite{ballard2008meet}, while some geometrical ideas behind it and their relations to string theory are discussed in the popularizing book by Yau,
\emph{The Shape of Inner Space: String Theory and the Geometry of the Universe's Hidden Dimensions} \cite{yau2010shape}.

\bibliographystyle{plain}
\bibliography{brouwer}{}

\begin{thebibliography}{10}

\bibitem{ballard2008meet}
Matthew~Robert Ballard.
\newblock Meet homological mirror symmetry.
\newblock {\em arXiv preprint arXiv:0801.2014}, 2008.

\bibitem{bott1988morse}
Raoul Bott.
\newblock Morse theory indomitable.
\newblock {\em Publications Math{\'e}matiques de l'IH{\'E}S}, 68:99--114, 1988.

\bibitem{brouwer1908continuous}
Luitzen Egbertus~Jan Brouwer.
\newblock On continuous vector distributions on surfaces.
\newblock {\em Koninklijke Nederlandse Akademie van Wetenschappen Proceedings
  Series B Physical Sciences}, 11:850--858, 1908.

\bibitem{brouwer1911abbildung}
Luitzen Egbertus~Jan Brouwer.
\newblock {\"U}ber abbildung von mannigfaltigkeiten.
\newblock {\em Mathematische Annalen}, 71(1):97--115, 1911.

\bibitem{cualduararu2005derived}
Andrei Caldd$\breve{\mathrm{a}}$r$\breve{\mathrm{a}}$ru.
\newblock Derived categories of sheaves: a skimming.
\newblock {\em Snowbird Lectures in Algebraic Geometry. Contemp. Math},
  388:43--75, 2005.

\bibitem{candelas1991pair}
Philip Candelas, C~De~La Ossa~Xenia, Paul~S Green, and Linda Parkes.
\newblock A pair of calabi-yau manifolds as an exactly soluble superconformal
  theory.
\newblock {\em Nuclear Physics B}, 359(1):21--74, 1991.

\bibitem{dieudonne2009history}
Jean Dieudonn{\'e}.
\newblock {\em A history of algebraic and differential topology, 1900-1960}.
\newblock Springer Science \& Business Media, 2009.

\bibitem{dieudonne1985history}
Suzanne~C Dieudonne.
\newblock {\em History Algebraic Geometry}.
\newblock CRC Press, 1985.

\bibitem{floer1986proof}
Andreas Floer.
\newblock Proof of the arnold conjecture for surfaces and generalizations to
  certain kaehler manifolds.
\newblock {\em Duke Mathematical Journal (C)}, 1986.

\bibitem{floer1988morse}
Andreas Floer.
\newblock Morse theory for lagrangian intersections.
\newblock {\em Journal of differential geometry}, 28(3):513--547, 1988.

\bibitem{freudenthal1976lej}
Hans Freudenthal.
\newblock Lej brouwer collected works.
\newblock 1976.

\bibitem{givental1996equivariant}
Alexander~B Givental.
\newblock Equivariant gromov-witten invariants.
\newblock {\em International Mathematics Research Notices}, 1996(13):613--663,
  1996.

\bibitem{guillemin2010differential}
Victor Guillemin and Alan Pollack.
\newblock {\em Differential topology}, volume 370.
\newblock American Mathematical Soc., 2010.

\bibitem{hofer2012floer}
Helmut Hofer, Clifford~H Taubes, Alan Weinstein, and Eduard Zehnder.
\newblock {\em The Floer memorial volume}, volume 133.
\newblock Birkh{\"a}user, 2012.

\bibitem{hutchings2002lecture}
Michael Hutchings.
\newblock Lecture notes on morse homology (with an eye towards floer theory and
  pseudoholomorphic curves).
\newblock {\em available at math. berkeley. edu/\~{}
  hutching/teach/276-2010/mfp. ps}, 2002.

\bibitem{kontsevich1994homological}
Maxim Kontsevich.
\newblock Homological algebra of mirror symmetry.
\newblock {\em arXiv preprint alg-geom/9411018}, 1994.

\bibitem{lian1998mirror}
B~Lian, Kefeng Liu, and Shing-Tung Yau.
\newblock Mirror principle, a survey.
\newblock {\em Current developments in mathematics (1998)(Cambridge, MA)},
  pages 35--82, 1998.

\bibitem{liu1998floer}
Gang Liu, Gang Tian, et~al.
\newblock Floer homology and arnold conjecture.
\newblock {\em J. Differential Geom}, 49(1):1--74, 1998.

\bibitem{miller2000leray}
Haynes Miller.
\newblock Leray in oflag xviia: the origins of sheaf theory, sheaf cohomology,
  and spectral sequences.
\newblock {\em Kantor 2000}, pages 17--34, 2000.

\bibitem{poincare2010papers}
Henri Poincar{\'e}.
\newblock {\em Papers on Topology: Analysis Situs and Its Five Supplements},
  volume~37.
\newblock American Mathematical Soc., 2010.

\bibitem{polishchuk2001categorical}
Alexander Polishchuk and Eric Zaslow.
\newblock Categorical mirror symmetry in the elliptic curve.
\newblock {\em AMS IP STUDIES IN ADVANCED MATHEMATICS}, 23:275--296, 2001.

\bibitem{salamon1999lectures}
Dietmar Salamon.
\newblock Lectures on floer homology.
\newblock {\em Symplectic geometry and topology (Park City, UT, 1997)},
  7:143--229, 1999.

\bibitem{smith2015symplectic}
Ivan Smith.
\newblock A symplectic prolegomenon.
\newblock {\em Bulletin of the American Mathematical Society}, 52(3):415--464,
  2015.

\bibitem{van2012lej}
Dirk Van~Dalen.
\newblock {\em LEJ Brouwer--Topologist, intuitionist, philosopher: How
  mathematics is rooted in life}.
\newblock Springer Science \& Business Media, 2012.

\bibitem{vonk2005mini}
Marcel Vonk.
\newblock A mini-course on topological strings.
\newblock {\em arXiv preprint hep-th/0504147}, 2005.

\bibitem{weibel1999history}
Charles~A Weibel.
\newblock {\em History of homological algebra}.
\newblock na, 1999.

\bibitem{witten1982supersymmetry}
Edward Witten et~al.
\newblock Supersymmetry and morse theory.
\newblock {\em J. diff. geom}, 17(4):661--692, 1982.

\bibitem{yau2010shape}
Shing-Tung Yau and Steve Nadis.
\newblock {\em The shape of inner space: string theory and the geometry of the
  Universe's hidden dimensions}.
\newblock Basic Books, 2010.

\end{thebibliography}

\end{document}